\documentclass[twoside,12pt,leqno]{article}
\advance\topmargin by -\headheight\advance\topmargin by -\headsep

\usepackage{amsmath}
\usepackage{amsthm}
\usepackage{amssymb}
\usepackage{amscd}
\usepackage{graphicx}
\usepackage{epsfig}

\numberwithin{equation}{section}
\allowdisplaybreaks

\newtheorem{theorem}{Theorem}
\newtheorem{lemma}{Lemma}[section]
\newtheorem{proposition}{Proposition}[section]
\newtheorem{corollary}{Corollary}[section]

\theoremstyle{definition}
\newtheorem{definition}{Definition}[section]
\newtheorem{example}{Example} [section]

\newtheorem{remark}{Remark}[section]

\usepackage{diagrams}
\usepackage{hyperref}

\def\co{\colon\thinspace}

\newcommand{\gequ}{\geqslant}

\DeclareMathOperator{\Cocyl}{Cocyl} \DeclareMathOperator{\Cyl}{Cyl}
\DeclareMathOperator{\Cocone}{Cocon} \DeclareMathOperator{\Cone}{Con}

\DeclareMathOperator{\End}{End} \DeclareMathOperator{\Der}{Der}
\DeclareMathOperator{\Vect}{Vect} \DeclareMathOperator{\ad}{ad}
\DeclareMathOperator{\Ker}{Ker}
\newcommand{\Imm}{{\mathop{\mathrm{Im}}}}
\newcommand{\void}{\varnothing}

\newcommand{\id}{{\mathrm{id}}}

\newcommand{\der}[2]{{\frac{\partial {#1}}{\partial {#2}}}}

\newcommand{\Z}{{\mathbb Z_{2}}}
\newcommand{\ZZ}{{\mathbb Z}}

\newcommand{\fun}{C^{\infty}}
\newcommand{\lsch}{{[\![}}
\newcommand{\rsch}{{]\!]}}

\def\e{\varepsilon}
\def\s{\sigma}

\def\D{\Delta}

\newcommand{\g}{{\gamma}}

\newcommand{\x}{{\xi}}

\newcommand{\LL}{{\Lambda}}

\newcommand{\at}{{\tilde a}}
\newcommand{\bt}{{\tilde b}}

\newcommand{\xt}{{\tilde x}}
\newcommand{\yt}{{\tilde y}}

\begin{document}

\setcounter{page}{1}
\thispagestyle{empty}



\markboth{Theodore Th. Voronov}{Higher Derived Brackets for
Arbitrary Derivations} 

\label{firstpage}
$ $
\bigskip

\bigskip

\centerline{{\Large   Higher Derived Brackets for Arbitrary
Derivations}}

\bigskip
\bigskip
\centerline{{\large by  Theodore Th. Voronov}}

\vspace*{.7cm}

\begin{abstract}
We introduce and study a construction of higher derived brackets
generated by a (not necessarily inner) derivation of a Lie
superalgebra. Higher derived brackets  generated by an element of a
Lie superalgebra were introduced  in our earlier work. Examples of
higher derived brackets naturally appear in geometry and mathematical
physics.  From a totally different viewpoint, we show that higher
derived brackets arise when one wants to turn the inclusion map of a
subalgebra of a differential Lie superalgebra, with a given
complementary subalgebra, into a fibration. (For a non-Abelian
complementary subalgebra, this leads to a generalization of
$L_{\infty}$-algebras with dropped or weakened (anti)symmetry of the
brackets.)
\end{abstract}

\pagestyle{myheadings}
\section{Introduction}

Higher derived brackets were introduced by the author
in~\cite{tv:higherder}, motivated by  physical and
differential-geometric examples. The starting point in the
construction was an element $\D$ in a Lie superalgebra $L$ provided
with a direct sum decomposition $L=K\oplus V$ into two subalgebras,
where $V$ is Abelian. Then a sequence of symmetric brackets on $V$
is `derived' from $\D$ in the same way as the partial derivatives of
a function:
\begin{equation*}
\{\underbrace{\x,\ldots,\x}_{n}\}_{\D}=P(-\ad \x)^n\D,
\end{equation*}
for coinciding even arguments, $P$ denoting the projector on $V$.
(In particular examples this analogy becomes exact.) It was proved
that the Jacobiators for the higher derived brackets of an odd $\D$
are equal to the higher derived brackets of $\D^2$. In particular,
this leads to $L_{\infty}$-algebras and algebras related with them.

In this paper we introduce and study higher derived brackets
generated by an arbitrary derivation $D\co L\to L$, which does not
have to be inner. See formula~\eqref{eqdefderbrack}. (The case of
non-inner derivations was touched on in the final version
of~\cite{tv:higherder}  without proofs.) As in~\cite{tv:higherder},
we make use of the decomposition $L=K\oplus V$. The subalgebra $V$ is
still assumed to be Abelian, though at the end of the paper we
briefly discuss how this condition can be relaxed.

Our first main result is Theorem~\ref{thmjac}, which we state and
prove in Section~\ref{secjac}. It says that the Jacobiators for the
higher derived  brackets generated by an odd derivation $D$ are
equal to the brackets generated by $D^2$. So it is an analog of a
similar statement for $\D$. However, the presently  available proof
of Theorem~\ref{thmjac} is technically much harder. Notice also that
strictly speaking, the theorem about brackets generated by $\D$ is
not a corollary of the theorem for $D$ because of the possible
presence of a `background' in the former.

Secondly, we establish a relation between higher derived brackets
and homotopical algebra. This is done in Section~\ref{sechomot}. The
main result is Theorem~\ref{thmcocyl}. The question of whether the
higher derived brackets of $\D$ defined in our
paper~\cite{tv:higherder} can be interpreted in the framework of
homotopical algebra was asked by the anonymous referee
of~\cite{tv:higherder}. In fact, he suggested linking them with the
notion of a `left cone' (i.e., a cocone, or a homotopy fiber in
topologists' language). This question happened to be very fruitful.
The proper framework for it is when the brackets are generated by an
arbitrary derivation $D$.  In Section~\ref{sechomot} we show that
such a homotopical-algebraic interpretation is indeed possible.  We
show that the higher derived brackets of $D$ appear as part of the
brackets  in $\Pi L\oplus V$ that naturally arise  from the
condition that the canonical differential in $\Pi L\oplus V$ (viewed
as a cone or a cocylinder) respects an algebra structure extended
from $L$, and we prove that the latter brackets make $\Pi L\oplus V$
an $L_{\infty}$-algebra  if $D^2=0$. Thus we arrive at an
alternative approach to  higher derived brackets.

Behind Theorem~\ref{thmjac} one can recognize a more general
algebraic statement. If one considers the Lie superalgebra $\Der L$
of derivations of $L$  and the Lie superalgebra $\Vect V$ of vector
fields on $V$, both w.r.t. the commutator, then it is possible to
see that the construction of higher derived brackets gives   a
homomorphism $\Der L\to \Vect V$. By identifying vector fields with
multilinear operations on $V$ specified by their Taylor expansion at
zero, we arrive at the statement that $V$ becomes a `generalized'
$L_{\infty}$-algebra `over' the Lie superalgebra $\Der L$ (that is,
there is a family of brackets parametrized by elements of $\Der L$
obeying `Jacobi type' relations following the relations in $\Der
L$). This is a new algebraic notion. We discuss it briefly. As
mentioned, we also briefly discuss the possibility of dropping the
condition that $V$ be Abelian. By doing so, we  arrive at higher
derived brackets that are not necessarily symmetric. This leads to
another generalization of $L_{\infty}$-algebras, which we hope to
analyze further elsewhere.

\textit{Terminology and notation.}  We use the `super' language and
conventions; in particular, a vector space always means a
$\Z$-graded vector space, and we freely identify it with the
corresponding vector supermanifold; multilinearity, symmetry,
antisymmetry, derivations, etc., are always understood in the
$\Z$-graded sense. $\Pi$ stands for the parity reversion functor,
and the parity of homogeneous objects is denoted by a tilde, i.e.
$\at=0$ or $\at=1$ if $a\in A_0$ or $a\in A_1$ respectively, for $a$
in a $\Z$-graded  module $A$.

\textit{Acknowledgements.} I wish to thank the anonymous referee of
paper~\cite{tv:higherder} for remarks that motivated the
``homotopical'' part of this work, and Hovhannes Khudaverdian for
numerous fruitful discussions. I am very grateful to Peter Symonds
as well as to the referee of the present paper for their help in
improving the style of exposition.

\section{Construction of Higher Derived Brackets}\label{secjac}

The algebraic setup is as follows. We are given a Lie superalgebra
$L$ and a decomposition of $L$ into a sum of two subalgebras:
\begin{equation*}
L=K\oplus V.
\end{equation*}
Let  $P\co L\to V$ be the projector on $V$ parallel to $K$, i.e.,
$V=\Imm P$, $K=\Ker P$.

Consider an arbitrary derivation $D\co L\to L$, either even or odd.

\begin{definition} \label{def.brack}
The {\em $k$-th (higher) derived bracket of $D$} is a multilinear
operation
$$
\underbrace{V\times \ldots \times V}_{\text{$k$ times}}\to V
$$
given by the formula
\begin{equation}\label{eqdefderbrack}
\{a_1,\ldots,a_k\}_D:=P[\ldots[[Da_1,a_2],a_3],\ldots,a_k]
\end{equation}
where $a_i\in V$. Here $k=1,2,3,\ldots \ $.
\end{definition}

The derived brackets have the same parity as the parity of the
derivation $D$.

In this paper we construct higher derived bracket for an arbitrary
derivation  $D$. Higher derived brackets  were first defined
in~\cite{tv:higherder} in the case when $D$ is an inner derivation,
$D=\ad \D$ for some $\D\in L$.

\begin{remark} The binary derived bracket when $P=\id$, i.e, $L=V$, $K=0$,
was introduced by Yvette~Kosmann-Schwarzbach~\cite{yvette:derived}
following an idea of Koszul, and independently by the
author~\cite{tv:lectures93etc} (unpublished).
In~\cite{yvette:derived} a slightly more general setting was
considered,  $L$ being a Loday (Leibniz) algebra. Derived brackets
have numerous applications, see~\cite{yvette:derived2} for a survey.
\end{remark}

\begin{lemma} Suppose that $V$ is an Abelian subalgebra.
Then the derived brackets  are symmetric (in  the $\Z$-graded sense).
\end{lemma}

\textbf{From now on we assume that $V$ is Abelian.}

(Later we shall discuss whether this requirement  can be relaxed.)

\begin{remark}
A symmetric multilinear operation is defined by its value on
coinciding even arguments (to be more precise, this is true if
extending scalars to include as many `odd constants' as necessary, is
allowed). For the higher derived brackets, if $\x$ is even,  we have
\begin{equation}\label{eqxxx}
\{\underbrace{\x,\ldots,\x}_{n}\}_D=(-1)^{n-1}P\left(\ad
\x\right)^{n-1} D\x
\end{equation}
for any $n=1,2,\ldots \ $, regardless of the parity of $D$.
\end{remark}

We want to investigate if in addition to symmetry the derived
brackets can satisfy other identities such as the Jacobi identity.
Consider the binary bracket. Notice first that it is symmetric, not
antisymmetric, compared to the bracket on a Lie algebra. To turn
symmetry into antisymmetry we have to reverse parity and consider
$\Pi V$.  The bracket induced on $\Pi V$ will be even (as a Lie
bracket should be) if the bracket in $V$ is odd.  Therefore, to ask
about (analogs of) the Jacobi identity for the higher derived
brackets makes sense when the derivation $D$ is odd.

\begin{example} \label{sec_const.ex1}
Consider a vector space $V$. Let $L=\Der S(V^*)$. Elements of $L$
can be viewed as polynomial vector fields on $V$. Let $\xi^i$ be the
linear coordinates on $V$ corresponding to a basis $(e_i)$ in $V$.
We can consider vectors in $V$ as vector fields with constant
coefficients, so $V\subset L$ will be an Abelian subalgebra. Then
$L=K\oplus V$, where the subalgebra $K$ consists of all vector
fields vanishing at the origin. The projector $P$ maps every vector
field to its value at the origin (constant vector field). Consider
an arbitrary  vector field
\begin{equation*}
X=X^j(\x)\,\der{}{\x^j},
\end{equation*}
even or odd, and consider the derived brackets on $V$ generated by
the derivation $\ad X$,
\begin{equation*}
\{u_1,\ldots,u_k\}_{X}:=\{u_1,\ldots,u_k\}_{\ad
X}=\left[\ldots\left[X,u_1\right],\ldots,u_k\right](0)\,.
\end{equation*}
One can see that
\begin{equation*}
\{e_{i_1},\ldots,e_{i_k}\}_{X}=(\pm)\, X^j_{i_1 \ldots i_k}\,e_j
\end{equation*}
where
\begin{equation*}
X^j_{i_1 \ldots i_k}=\frac{\partial^k X^j}{\partial\x^{i_1}\ldots
\partial\x^{i_k}}\,(0)\,.
\end{equation*}
In particular, consider a quadratic vector field:
\begin{equation*}
Q=\frac{1}{2}\,\xi^i\xi^j Q_{ji}^k\der{}{\x^k}\,.
\end{equation*}
Then the $k$-th derived bracket of $Q$ is zero unless $k=2$, and the
binary derived bracket is given by
\begin{equation}\label{eqbrackqbin}
\{e_i,e_j\}_Q=(\pm)\,Q_{ij}^ke_k\,.
\end{equation}
Suppose that $Q$ is odd. By a direct check one can see that
\textit{the Jacobi identity for the bracket~\eqref{eqbrackqbin} is
equivalent to the condition $Q^2=0$.} (More precisely, the usual
graded Jacobi identity is valid for the antisymmetric bracket in
$\Pi V$.) $Q$ then can be identified with the Chevalley--Eilenberg
differential in the standard cochain complex of the resulting Lie
(super)algebra.
\end{example}

Odd vector fields with square zero are known as
\textit{homological}. Example~\ref{sec_const.ex1} shows that the
structure of a Lie superalgebra on a vector space corresponds  to a
quadratic homological vector field. If we drop the condition that
the homological vector field be quadratic, we obtain
`$L_{\infty}$-algebras' or `strong homotopy Lie algebras' where the
Jacobi identity for a binary bracket holds up to homotopy, which is
a ternary bracket, and in its turn satisfies an analog of the Jacobi
identity up to a homotopy, and so on. Describing algebraic
structures by derivations of square zero is a very general principle
dating back to Nijenhuis in 1950's and  used in recent works of
Kontsevich (his ``formality theorem'' implying the existence of
deformation quantization of Poisson manifolds,
see~\cite{kontsevich:quant}).

Let $V$ be a vector space endowed with a sequence of $k$-linear odd
symmetric operations denoted by braces. Here $k=0,1,2,3,\ldots \ $.
\begin{definition} \label{def.jacobiator} The \textit{$n$-th Jacobiator} is the following expression:
\begin{equation*}
J^n(a_1,\ldots,a_n)=\sum_{k+l=n}\,\sum_{\text{$(k,l)$-shuffles}}
(-1)^{\e}
\{\{a_{\s(1)},\ldots,a_{\s(k)}\},a_{\s(k+1)},\ldots,a_{\s(k+l)}\}
\end{equation*}
where the sign $(-1)^{\e}$ is given by the usual sign rule for
permutations of homogeneous elements of $V$.
\end{definition}

Let us recall the  definition of $L_{\infty}$-algebras (due to Lada
and Stasheff~\cite{lada:stasheff}).
\begin{definition}[Stasheff, Lada]
An \textit{$L_{\infty}$-algebra}, or \textit{strong(ly) homotopy Lie
algebra}, is   a vector space $V$  endowed with a sequence of
$k$-linear odd symmetric operations, $k=0,1,2,3,\ldots \ $, such that
all the Jacobiators vanish.
\end{definition}

\begin{remark} We gave the definition in the form most convenient for our purposes.
If one wishes to directly include the case of ordinary Lie algebras,
the term $L_{\infty}$-algebra should be applied to the structure on
the opposite space, i.e., $\Pi V$, where the corresponding
operations are antisymmetric and are even for an even number of
arguments and odd otherwise. Also, it is often assumed that the
$0$-bracket is zero. The $0$-bracket is sometimes referred to as the
`background'.
\end{remark}

\begin{proposition} \label{prop.linfandhom}There is a one-to-one correspondence between
$L_{\infty}$-algebra structures on $V$ and formal homological vector
fields on $V$:
\begin{equation*}
Q=Q^k(\xi)\,\der{}{\x^k}=\left(Q^k_0+\x^iQ_i^k+\frac{1}{2}\,\xi^i\xi^j
Q_{ji}^k+\frac{1}{3!}\,\xi^i\xi^j\x^l Q_{lji}^k+\ldots
\right)\der{}{\x^k}\,.
\end{equation*}
\end{proposition}

This proposition is a well-known fact. What we can give here is an
explicit invariant expression for the  correspondence: \textit{the
brackets in an $L_{\infty}$-algebra corresponding to a homological
field $Q$ are the higher derived brackets of $Q$},
\begin{equation}\label{equuQ}
\{u_1,\ldots,u_k\}_{Q}=\left[\ldots\left[Q,u_1\right],\ldots,u_k\right](0)
\end{equation}
(this generalizes~Example~\ref{sec_const.ex1}).

Let us return to our abstract setup. Consider the higher derived
brackets~\eqref{eqdefderbrack} of an \textit{odd} derivation $D$. We
get a sequence of odd symmetric multilinear operations on $V$. By
definition  the $0$-ary operation is zero. What about the
Jacobiators?

\begin{theorem}
\label{thmjac} Suppose that $D$ preserves the subalgebra $K=\Ker P$.
Then the $n$-th Jacobiator of the derived brackets of $D$ equals the
$n$-th derived bracket of $D^2$:
\begin{equation} \label{eqjacobimain}
J^n_D(a_1,\ldots,a_n)=\{a_1,\ldots,a_n\}_{D^2}\,,
\end{equation}
for all $n=1,2,3,\ldots \ $.
\end{theorem}

Let us make two comments before giving the proof.

Firstly, notice that in our setup the condition that $L=K\oplus V$ is
the sum of subalgebras where $V$ is Abelian can be expressed by the
identities
\begin{equation}\label{eqvabel}
    [Pa, Pb]=0
\end{equation}
and
\begin{equation}\label{eqdistrib}
    P[a,b]=P[Pa,b]+P[a,Pb],
\end{equation}
for all $a,b$ (a ``distributivity law'' for $P$). Notice
that~\eqref{eqdistrib} is also equivalent to vanishing of the
Nijenhuis bracket of $P$ with itself.

Secondly,  the condition in the theorem that $D(K)\subset K$ can be
written as the identity
\begin{equation}\label{eqpdp}
PDP=PD\,.
\end{equation}
Condition~\eqref{eqpdp}  already appears if we check the  Jacobiator
of order one:
\begin{equation*}
J^1_D(a)=\{\{a\}_D\}_D=PDPDa=PDDa=\{a\}_{D^2},
\end{equation*}
if $PDP=PD$. (Notice that  in general $D$  does not have to preserve
$V$. Indeed, if $D$ preserves $V$, then all the derived brackets of
$D$ starting from the binary bracket, will vanish.)

\begin{proof}[Proof of Theorem~\ref{thmjac}]
To simplify the notation, let us omit temporarily the subscript $D$
from the brackets and Jacobiators.  Since the Jacobiators are
multilinear symmetric functions, it is sufficient to consider them
for coinciding even arguments. Denote $J^n(\x, \ldots, \x)$ where
$\x$ is even by $J^n(\x)$. From Definition~\ref{def.jacobiator} and
equation~\eqref{eqxxx} we clearly obtain
\begin{equation*}
J^n(\x)= \sum_{l=0}^{n-1}C_n^l
\{\{\underbrace{\x,\ldots,\x}_{n-l}\},\underbrace{\x,\ldots,\x}_{l}\},
\end{equation*}
where $C_n^l=\frac{n!}{l!(n-l)!}$ is the binomial coefficient, in our
case appearing as the number of $(n-l,l)$-shuffles. It follows that
\begin{multline*}
    J^n(\x)=
\sum_{l=0}^{n-1}C_n^lP[\ldots[D\{\underbrace{\x,\ldots,\x}_{n-l}\},
\underbrace{\x],\ldots,\x]}_{l}=\\
\sum_{l=0}^{n-1}C_n^lP[\ldots[DP(-1)^{n-l-1}
(\ad\x)^{n-l-1}D\x,\underbrace{\x],\ldots,\x]}_{l}=\\
\sum_{l=0}^{n-1}C_n^l (-1)^{n-l-1} (-1)^l P(\ad \x)^l
DP(\ad\x)^{n-l-1}D\x=\\
\sum_{l=0}^{n-1}C_n^l (-1)^{n-1}  P(\ad \x)^l DP(\ad\x)^{n-l-1}D\x\,.
\end{multline*}
Consider $(\ad \x)^l DP$. We want to move $D$ to the left. Since
$$
\ad\x\cdot D-D\cdot \ad\x=-\ad(D\x)\,,
$$
as one can easily check, it follows that for any $l\gequ 1$
\begin{multline*}
    (\ad \x)^l DP=(\ad\x)^{l-1}(\ad\x\cdot D-D\cdot \ad\x+D\cdot
    \ad\x)P=\\
    (\ad\x)^{l-1}(-\ad(D\x)+D\cdot
    \ad\x)P=-(\ad\x)^{l-1}\ad(D\x)P=
    -\left[(\ad\x)^{l-1}D\x,  P(.)
    \right]
\end{multline*}
where we used $\ad\x\cdot P=0$. Substituting this into the formula
above we obtain
\begin{equation*}
    J^n(\x)=\sum_{l=1}^{n-1} C_n^l(-1)^nP\left[(\ad\x)^{l-1}D\x,
    P(\ad\x)^{n-l-1}D\x\right]
+(-1)^{n-1}PDP(\ad\x)^{n-1}D\x
\end{equation*}
or
\begin{equation*}
   (-1)^n J^n(\x)=\sum_{l=1}^{n-1} C_n^l P\left[(\ad\x)^{l-1}D\x,
    P(\ad\x)^{n-l-1}D\x\right]
-PD(\ad\x)^{n-1}D\x
\end{equation*}
(where we also used~\eqref{eqpdp}). We can re-arrange the first sum
by adding it to itself in the reverse order and dividing by two:
\begin{multline*}
   \sum_{l=1}^{n-1} C_n^l P\left[(\ad\x)^{l-1}D\x,
    P(\ad\x)^{n-l-1}D\x\right]=\\
    \frac{1}{2}\,\left(\sum_{l=1}^{n-1} C_n^l P\left[(\ad\x)^{l-1}D\x,
    P(\ad\x)^{n-l-1}D\x\right]+
    \sum_{l=1}^{n-1} C_n^l P\left[(\ad\x)^{n-l-1}D\x,
    P(\ad\x)^{l-1}D\x\right]\right)=\\
    \frac{1}{2}\,\sum_{l=1}^{n-1} C_n^l  \Bigl(P\left[(\ad\x)^{l-1}D\x,
    P(\ad\x)^{n-l-1}D\x\right]+ P\left[(\ad\x)^{n-l-1}D\x,
    P(\ad\x)^{l-1}D\x\right]\Bigr)\,.
\end{multline*}
Noticing that $\left[(\ad\x)^{n-l-1}D\x,P(\ad\x)^{l-1}D\x\right]=
\left[P(\ad\x)^{l-1}D\x,(\ad\x)^{n-l-1}D\x\right]$, because $\x$ is
even and $D\x$ is odd, and using the distributivity
relation~\eqref{eqdistrib}, we find the following expression for the
Jacobiator:
\begin{equation*}
   (-1)^n J^n(\x)=
   \frac{1}{2}\,\sum_{l=1}^{n-1} C_n^l  P\left[(\ad\x)^{l-1}D\x,
    (\ad\x)^{n-l-1}D\x\right]-PD(\ad\x)^{n-1}D\x
\end{equation*}
or
\begin{multline}\label{eqjacob1}
   (-1)^n J^n(\x)=
\frac{1}{2}\,P\sum_{l=1}^{n-1} C_n^l  \left[(\ad\x)^{l-1}D\x,
    (\ad\x)^{n-l-1}D\x\right]-\\
    P\left[D,(\ad\x)^{n-1}\right]D\x- P(\ad\x)^{n-1}D^2\x\,.
\end{multline}
We shall now analyze the term $\left[D,(\ad\x)^{n-1}\right]D\x$.
Using the formula for the  commutator $[A,B^N]$ for arbitrary
operators $A$, $B$, we get
\begin{multline*}
    \left[D,(\ad\x)^{n-1}\right]D\x=\sum_{i+j=n-2}(\ad\x)^i[D,\ad\x](\ad\x)^jD\x= \\
    \sum_{i+j=n-2}(\ad\x)^i\ad(D\x)(\ad\x)^jD\x=
    \sum_{i+j=n-2}(\ad\x)^i\left[D\x, (\ad\x)^jD\x\right]=\\
    \sum_{i+j=n-2}\,\,\sum_{r+s=i}
    C_i^r\left[(\ad\x)^{r}D\x, (\ad\x)^{s+j}D\x\right]=
\sum_{i=0}^{n-2}\sum_{r=0}^i
    C_i^r\left[(\ad\x)^{r}D\x, (\ad\x)^{n-2-r}D\x\right]=\\
\sum_{r=0}^{n-2}\sum_{i=r}^{n-2}
    C_i^r\left[(\ad\x)^{r}D\x, (\ad\x)^{n-2-r}D\x\right]\,.
\end{multline*}
Since in the internal sum the commutators do not depend on the index
of summation  $i$, they can be taken out of the sum. It is possible
to apply a well-known identity for  sums of binomial coefficients
(see, e.g.~\cite[p.~153]{bron-sem:1981}):
\begin{equation*}
\sum_{i=r}^{m} C_i^r=C_r^r+C_{r+1}^r+\ldots+C_{m}^r=
C_r^0+C_{r+1}^1+\ldots+C_{m}^{m-r}=C_{m+1}^{m-r},
\end{equation*}
where in our case $m=n-2$. Hence
\begin{equation*}
\sum_{i=r}^{n-2} C_i^r=C_{n-1}^{n-2-r},
\end{equation*}
and we arrive at the equality
\begin{equation}\label{eqkomdiad1}
    \left[D,(\ad\x)^{n-1}\right]D\x=
\sum_{r=0}^{n-2}C_{n-1}^{n-2-r}\left[(\ad\x)^{r}D\x,
(\ad\x)^{n-2-r}D\x\right]\,.
\end{equation}
Notice that since $D\x$ is odd and the bracket is symmetric, the
left-hand side contains similar terms, with $r$ and $r'$, where
$r=n-2-r'$. Hence, this sum can be re-arranged by writing it twice in
opposite orders and dividing by two:
\begin{multline*}
    \left[D,(\ad\x)^{n-1}\right]D\x= \\
\frac{1}{2}\left(\sum_{r=0}^{n-2}C_{n-1}^{n-2-r}\left[(\ad\x)^{r}D\x,
(\ad\x)^{n-2-r}D\x\right]+
\sum_{r=0}^{n-2}C_{n-1}^{r}\left[(\ad\x)^{n-2-r}D\x,
(\ad\x)^{r}D\x\right]\right)=\\
\frac{1}{2}\sum_{r=0}^{n-2}\left(C_{n-1}^{r+1}+C_{n-1}^{r}\right)
\left[(\ad\x)^{r}D\x,(\ad\x)^{n-2-r}D\x\right]=\\
\frac{1}{2}\sum_{r=0}^{n-2} C_{n}^{r+1}
\left[(\ad\x)^{r}D\x,(\ad\x)^{n-2-r}D\x\right]=
\frac{1}{2}\sum_{l=1}^{n-1} C_{n}^{l}
\left[(\ad\x)^{l-1}D\x,(\ad\x)^{n-1-l}D\x\right]\,,
\end{multline*}
which coincides with the first term in the formula for the
Jacobiator~\eqref{eqjacob1}. Substituting into~\eqref{eqjacob1}, we
see that the first two terms cancel, and we finally obtain
\begin{equation*}
(-1)^n J^n(\x)=-P(\ad\x)^{n-1}D^2\x
\end{equation*}
or
\begin{equation*}
J^n_D(\x)=(-1)^{n-1}
P(\ad\x)^{n-1}D^2\x=\{\underbrace{\x,\ldots,\x}_{n}\}_{D^2}
\end{equation*}
for an arbitrary even $\x$. This implies
identity~\eqref{eqjacobimain} for all elements of $V$.
\end{proof}

We say that the derivation $D$ is of \textit{order} $r$ with respect
to the subalgebra $V$ if for all $a_1,\ldots, a_{r+1}\in V$
\begin{equation*}
\left[\ldots\left[Da_1,a_2\right],\ldots,a_{r+1}\right]=0.
\end{equation*}
Here $r=0,1,2,\ldots\,$

If $D$ is of order $r$, all the derived $k$-brackets of $D$ vanish
for $k\gequ r+1$.

\begin{corollary} For an odd derivation $D$, if the order of $D^2$ is $r$, then
the higher derived brackets of  $D$ satisfy all the Jacobi
identities with $\gequ r+1$ arguments,
\end{corollary}

\begin{corollary} If the order of $D^2$ is zero, i.e., $D^2(V)=0$, in particular if
$D^2=0$, then the higher derived brackets of an odd derivation $D$
define an $L_{\infty}$-algebra.
\end{corollary}

Proposition~\ref{prop.linfandhom} shows that all
$L_{\infty}$-algebras are obtained in this way.

\section{Examples}

All examples of higher derived brackets naturally arising in
applications are for the case when $D=\ad \D$ is an inner derivation
given by some element $\D$. This is the situation where higher
derived brackets were first introduced in~\cite{tv:higherder}. An
analog of Theorem~\ref{thmjac} was proved there for brackets
generated by $\D$. (That proof is simpler than the above proof of
Theorem~\ref{thmjac} for general $D$.) Let us clarify the relation
between the higher derived brackets of an element $\D\in L$ as
introduced in~\cite{tv:higherder} and the higher derived brackets of
a derivation $D\co L\to L$ as in Definition~\ref{def.brack}.

Any element $\D\in L$, of course, gives an inner derivation $D=\ad
\D\co L\to L$, and the higher derived brackets of the derivation
$D=\ad \D$
\begin{equation*}
\{a_1,\ldots,a_k\}_{D}=P\left[\ldots\left[(\ad
\D)a_1,a_2\right],\ldots,a_k\right],
\end{equation*}
coincide with the brackets defined in~\cite{tv:higherder},
\begin{equation*}
\{a_1,\ldots,a_k\}_{\D}=P\left[\ldots\left[\left[\D,a_1\right],a_2\right],
\ldots,a_k\right],
\end{equation*}
where $k=1,2,3,\ldots \ $. However, for $\D$ there is a natural
notion of a $0$-bracket (no arguments, a distinguished element),
\begin{equation*}
\{\void\}_{\D}=P\D\,,
\end{equation*}
which is not defined for  arbitrary derivations $D$. The Jacobiators
for the higher derived brackets of $\D$ include this $0$-bracket and
start with the $0$-th Jacobiator $\{\{\void\}_{\D}\}_{\D}$. At the
same time, the $0$-ary bracket  is assumed to be zero in all the
Jacobiators for a general $D$ and it does not appear in our
Theorem~\ref{thmjac}. There is no obvious way of incorporating the
$0$-th bracket into the picture for a general derivation $D$. If
$P\D\neq 0$, that means that $\D\notin K$, hence there is no
guarantee that $\ad\D(K)\subset K$, which is a condition of
Theorem~\ref{thmjac}. The calculation of $J^1_D(a)$ above shows that
some sort of condition is needed (and at least a weaker condition
$PD^2P=PDPDP$ is necessary). Therefore, Theorem~1
of~\cite{tv:higherder}, to which  Theorem~\ref{thmjac} is an analog,
does not follow from Theorem~\ref{thmjac}, in general.

We can summarize by saying that the theory developed
in~\cite{tv:higherder} is a particular case of the theory developed
here if $P(\D)=0$, i.e., $\D\in K$. Then, in particular, $(\ad
\D)(K)\subset K$ and Theorem~\ref{thmjac} applies.

We shall leave open the question of how the theory for non-inner
derivations can be modified to incorporate an $0$-ary bracket.

With having this in mind, there are some examples of higher derived
brackets, all coming from inner derivations. They are given for
illustrative purposes only. More details can be found
in~\cite{tv:higherder}. See also~\cite{tv:graded},
\cite{tv:laplace2}.

\begin{example} The setup of Example~\ref{sec_const.ex1}.
$L=\Vect V$, where $V$ is a vector space, $P\co X\mapsto X(0)$ is a
projection onto the Abelian subalgebra of vector fields with
constant coefficients. For an odd vector field $Q$ such that
$Q(0)=0$ we get the higher derived brackets on $V$, $k=1,2,\ldots, \
$,
\begin{equation*}
\{u_1,\ldots,u_k\}_Q=\left[\ldots\left[Q,u_1\right],\ldots,u_k\right](0).
\end{equation*}
They define an $L_{\infty}$-algebra with e zero background (`strict'
in the terminology of~\cite{tv:higherder}) if $Q^2=0$, and this is a
canonical description of all (strict) $L_{\infty}$-algebra
structures on the space $V$.
\end{example}

\begin{example} $L=\End A$ for a commutative associative algebra with
unit $A$ and $V=A$. The projector $P$ maps an operator $\D$ to $\D
1\in A\subset \End A$. The higher derived brackets of $\ad\D$ for an
operator $\D$ such that $\D1=0$ are the `Koszul operations'
(see~\cite{koszul:crochet85})
\begin{equation*}
\{a_1,\ldots,a_k\}_{\D}=\left[\ldots\left[\D,a_1\right],\ldots,a_k\right]1,
\end{equation*}
$k=1,2,3,\ldots, ...\ $. For a differential operator of order $n$ the
brackets with more than $n$ arguments vanish and the top bracket is
the symbol of $\D$. An odd operator $\D$ satisfying $\D^2=0$ provides
an example of a  `homotopy Batalin--Vilkovisky algebra'.
\end{example}

\begin{example}
$L=\fun(T^*M)$, $V=\fun(M)$,  $P$ is the restriction on $M$, and
$i^*\co\fun(T^*M)\to\fun(M)$, where $i\co M\to T^*M$. For a
Hamiltonian $S\in\fun(T^*M)$ such that $i^*S=0$, on functions on $M$
there are derived brackets
\begin{equation*}
\{f_1,\ldots,f_k\}_{S}=i^*\left(\ldots\left(S,f_1\right),\ldots,f_k\right),
\end{equation*}
$k=1,2,3,\ldots\ $, where in the right-hand side stand the canonical
Poisson brackets on $T^*M$. If $S$ is odd (for a nontrivial example
$M$ should be a supermanifold) and satisfies $(S,S)=0$, we get
`higher Schouten brackets' on $\fun(M)$ giving an example of a
`homotopy Schouten algebra'.
\end{example}

\begin{example} Similarly, let $L=\fun(\Pi T^*M)$, $V=\fun(M)$ and let $P$ be
the restriction on $M$. For a multivector field $\psi\in\fun(\Pi
T^*M)$ such that $i^*\psi =0$, on functions on $M$ there are derived
brackets
\begin{equation*}
\{f_1,\ldots,f_k\}_{\psi}=i^* \lsch\ldots \lsch \psi,f_1
\rsch,\ldots,f_k \rsch,
\end{equation*}
$k=1,2,3,\ldots \, $, where on the right-hand side we have the
canonical Schouten brackets on $\Pi T^*M$. Since the canonical
Schouten brackets are odd, for an even $\psi$ the derived brackets
have alternating parity (even for an even number of arguments, odd
for odd). If $\lsch \psi,\psi\rsch=0$, these `higher Poisson
brackets' on functions on $M$  give an example of a `homotopy
Poisson algebra'.
\end{example}

Other examples of higher derived brackets  which we shall not
consider here, are   `homotopy Lie algebroids', which are an analog
of $L_{\infty}$-algebras in the world of algebroids, and the
non-linear analogs in the world of graded manifolds~\cite{tv:graded}
(see also~\cite{roytenberg:graded}). We hope to be able to say more
about such examples elsewhere.

It is a good question whether a genuinely non-inner derivation can
naturally occur in examples of higher derived brackets coming from
differential geometry or physics.

\section{Relation with Homotopy Theory}\label{sechomot}

Now we shall show how our construction of the (higher) derived
brackets arises naturally if one wishes to consider the homotopy
theory of Lie superalgebras.

Let us re-formulate the setup in a way convenient for this purpose.
We have a Lie superalgebra $L$ with a decomposition $L=K\oplus V$,
where $K$ and $V$ are subalgebras. Consider an odd derivation $D$
such that $D(K)\subset K$, and from the start assume that $D$ is of
square zero. Hence we have an inclusion of differential Lie
superalgebras
\begin{equation*}
i\co K\to L
\end{equation*}
and a given complement for the image of $i$, which is called $V$.
($V$ is \textit{not}, in general, a differential subalgebra.)

There is an idea, familiar to topologists, that every map can be
made into a fibration by appropriately replacing a space by a
homotopy equivalent one. More precisely, if we have a category where
a ``weak equivalence'', ``fibration'' and ``cofibration'' make sense
(i.e., a Quillen model category~\cite{quillen:hoalgebra67}), the
following diagram is called a \textit{cocylinder diagram}:
\begin{equation*}
\begin{diagram}[small]
    X &       & \rTo^f          &       & Y  \\
      & \rdTo_j &               & \ruTo_p &    \\
      &       & Z &       &    \\
\end{diagram}
\end{equation*}
if $j$ is a cofibration and weak equivalence, and $p$ is a
fibration. Then $Z$ is also denoted by $\Cocyl f$. (To refresh the
intuition,  recall that for topological spaces that are cell
complexes, cofibrations are just inclusions of subspaces, fibrations
are `Serre fibrations', i.e., maps satisfying the covering homotopy
property, and weak equivalences are maps inducing isomorphisms of
all homotopy groups. In this case, a cocylinder for any map $f\co
X\to Y$ may be obtained as a subspace in $X\times Y^I$ consisting of
all pairs $(x,\g)$ where $\g\co I\to Y$ is a path such that
$\g(0)=f(x)$.)

Can we do this (in an algebraic context) for the inclusion $K\to L$?

To begin with let us temporarily forget about the algebra structure.
Consider just an arbitrary inclusion of complexes
\begin{equation*}
i\co K\to L
\end{equation*}
such that there is a given complementary subspace $V$ (not a
subcomplex!) and $L=K\oplus V$. In the context of this paper, a
\textit{complex} is simply a vector space with an odd operator of
square zero.

Let $P$ be the projector onto $V$ parallel to $K$. The space $V$
becomes a complex with the differential $PD$. Introduce into $L\oplus
\Pi V$ an operator $d$ as follows:
\begin{equation}\label{eqdincocyl}
    d(x,\Pi a):=\bigl(Dx, -\Pi P(x+Da)\bigr), 
\end{equation}
for $x\in L$, $a\in V$. It is then straightforward to show that
$d^2=0$. Consider the maps $j\co K\to L\oplus \Pi V$ and $p\co
L\oplus \Pi V\to L$, where $j\co x\mapsto (x,0)$, $p\co (x,\Pi
a)\mapsto x$.

\begin{lemma}\label{lemcocyl}
The following diagram
\begin{equation*}
\begin{diagram}[small]
    K &       & \rTo^i          &       & L  \\
      & \rdTo_j &               & \ruTo_p &    \\
      &       & L\oplus \Pi V &       &    \\
\end{diagram}
\end{equation*}
is a cocylinder diagram in the category of complexes, i.e., the maps
$j$ and $p$ are chain maps,   $i=p\circ j$, the map  $j\co K\to
L\oplus \Pi V$ is a monomorphism \textnormal{(`cofibration')} and a
quasi-isomorphism \textnormal{(`weak homotopy equivalence')}, and the
map $p\co L\oplus \Pi V\to L$ is an epimorphism
\textnormal{(`fibration')}.
\end{lemma}
\begin{proof} This is immediate. A quasi-inverse for $j$ is the map
$$
q\co (x,\Pi a)\mapsto (1-P)(x+Da).
$$
\end{proof}

\begin{remark} As   is well known, for maps of complexes there are canonical
constructions of  cylinders and cocylinders; they are modelled on
the (co)chain complexes corresponding to the canonical topological
cylinders and cocylinders.  For a particular chain map it might be
more convenient to consider a `smaller' cylinder or cocylinder than
the one featured by the standard construction. This is exactly what
happens in our case. The standard cocylinder construction applied to
the inclusion $i\co K\to L$ would  not yield the complex $L\oplus
\Pi V$ as in Lemma~\ref{lemcocyl}, instead it would give a bigger
complex $K\oplus L\oplus \Pi L=K\oplus K\oplus V\oplus \Pi K\oplus
\Pi V$ that is homotopy equivalent to $L\oplus \Pi V$. One should
also note that the complex $L\oplus \Pi V$ essentially coincides
with the standard (co)cone of the projection $L\to V$. See the
Appendix.
\end{remark}

It follows from Lemma~~\ref{lemcocyl} that the space $\Pi V=\Ker p$,
taken with the differential $-\Pi PD$, is a homotopy fiber of the
inclusion of complexes $i\co K\to L=K\oplus V$.

Now we want to `turn the algebra structure on'. To this end, since
we have been working  with $V$ rather than $\Pi V$, let us first
apply a parity shift to the cocylinder diagram above. Then we have
the cocylinder diagram
\begin{equation*}
\begin{diagram}[small]
    \Pi K &       & \rTo^i          &       & \Pi L  \\
      & \rdTo_j &               & \ruTo_p &    \\
      &       & \Pi L\oplus  V &       &    \\
\end{diagram}
\end{equation*}
for the inclusion of complexes $\Pi K\to \Pi L$.  In particular, the
differential in $\Pi L\oplus V$ is
\begin{equation}\label{eqdifincoc2}
    d\co (\Pi x,a)\mapsto (-\Pi Dx, P(x+Da))
\end{equation}
(which is the differential in the standard cone, see the Appendix, of
the projection of complexes $(L,D)$ onto $(V, PD)$).

Let us restore our framework. Assume as above that $L$ is a Lie
superalgebra with $D$ being a derivation, and that $V$ is an Abelian
subalgebra. The Lie bracket in $L$ induces an odd bracket in $\Pi
L$:
\begin{equation}\label{eqbrackpixpiy}
\{\Pi x, \Pi y\}=\Pi[x,y](-1)^{\tilde x}.
\end{equation}
Is it possible to extend this to a bracket on the whole of $\Pi
L\oplus V$?

\begin{proposition} There exist an  odd binary bracket on $\Pi L\oplus V$ extending
that on $\Pi L$ such that the operator~\eqref{eqdifincoc2} acts as a
derivation. It is given by the formulae
\begin{align}
 \{\Pi x, a\}&=P[x,a], \label{eqbrackpixa}\\
 \{a, b\}&= P[Da,b]\label{eqbrackab}
\end{align}
for arbitrary $x\in L$ and $a,b\in V$.
\end{proposition}
\begin{proof}
As a starting point we use  formula~\eqref{eqbrackpixpiy} for the
bracket on $\Pi L$, where $\Pi x$ and $\Pi y$ are considered as
elements of $\Pi L\oplus V$. Apply $d$ given by~\eqref{eqdifincoc2}
to $\{\Pi x,\Pi y\}$ and require that the Leibniz rule be satisfied:
\begin{equation}\label{eqleibpixpiy}
d\{\Pi x,\Pi y\}=-\{d\Pi x,\Pi y\}-(-1)^{\xt+1}\{\Pi x,d\Pi y\}
\end{equation}
for all $x,y\in L$ (notice that the parity in~\eqref{eqbrackpixpiy}
`sits' at the opening bracket, hence the signs). Expanding $d$
by~\eqref{eqdifincoc2}, so that $d\Pi x=-\Pi Dx+Px$, and treating
the brackets between elements of $\Pi L$ and $V$ as  unknown, we
find that the failure of~\eqref{eqleibpixpiy}  for $x=y$ and $\xt=1$
is the difference $\{Px,\Pi x\}-P[Px,x]$. Replacing $Px$ by an
arbitrary element of $V$, we  arrive at the above
definition~\eqref{eqbrackpixa}. Now assume~\eqref{eqbrackpixa} and
require  the Leibniz rule for this new bracket:
\begin{equation}\label{eqleibpixa}
    d\{\Pi x,b\}=-\{d\Pi x,b\}+(-1)^{\xt}\{\Pi x,db\}
\end{equation}
for all $x\in L$, $b\in V$. Here $d\Pi x=-\Pi Dx+Px$, $da=PDa$, and
we treat the bracket in $V$ as unknown. The failure
of~\eqref{eqleibpixa} equals
$\{Px,b\}+(-1)^{\xt}P\left[Px,Db\right]$. Denoting $Px=a\in V$, we
arrive at the formula $\{a,b\}=-(-1)^{\at}P[a,Db]$ or, equivalently,
\begin{equation*}
\{a,b\}=P[Da,b]
\end{equation*}
as the necessary and sufficient condition of~\eqref{eqleibpixa}.
This is our derived bracket~\eqref{eqdefderbrack} for $k=2$. The
Leibniz rule for $\{a,b\}$ is now satisfied automatically and does
not bring any new relations.
\end{proof}

\begin{remark}
Defining the bracket by formula~\eqref{eqbrackpixa} is a sufficient
condition for~\eqref{eqleibpixpiy}. A more detailed analysis shows
that~\eqref{eqbrackpixa} is also necessary at least when $x\in K$.
Hence the condition that the operator~\eqref{eqdifincoc2} acts as a
derivation defines the bracket in an essentially unique way.
\end{remark}

One can see that a binary bracket defined in this way on $\Pi
L\oplus V$ will not satisfy the Jacobi identity exactly, thus giving
rise to a ternary bracket, and so on. Define  the higher brackets on
$\Pi L\oplus V$ as follows:
\begin{align}
\{\Pi x, a_1,\ldots,a_n\}&=P[\ldots[x,a_1],\ldots,a_n], \label{eqbrackxaaa}\\
\{a_1,\ldots,a_n\}&=P[\ldots[Da_1,a_2],\ldots,a_n],\label{eqbrackaaa}
\end{align}
where $n\gequ 1$.  As an unary bracket   take the
differential~\eqref{eqdifincoc2},  and set the $0$-ary bracket to
zero.  All the other brackets except those obtainable  by symmetry,
are defined to be  zero. Formulae~\eqref{eqbrackxaaa},
\eqref{eqbrackaaa} directly extend~\eqref{eqbrackpixa},
\eqref{eqbrackab} to many arguments, and formula~\eqref{eqbrackaaa}
is the familiar higher derived bracket on $V$ for all $k$.

\begin{theorem}
\label{thmcocyl} The set of brackets~\eqref{eqbrackpixpiy},
\eqref{eqbrackxaaa} and \eqref{eqbrackaaa}, together
with~\eqref{eqdifincoc2}, define on the space $\Pi L\oplus V$ the
structure of an $L_{\infty}$-algebra.
\end{theorem}

\begin{proof}   We shall prove that all the
brackets~\eqref{eqdifincoc2}--\eqref{eqbrackaaa} satisfy all the
generalized Jacobi identities. Consider the Jacobiator $J^n$ in $\Pi
L\oplus V$ with $n$ arbitrary  arguments. Without loss of generality
we can assume that each of the arguments is either in $\Pi L$ or
$V$. We claim that there can be no non-trivial Jacobiators with more
than $3$ arguments in $\Pi L$. Indeed, $J^n$ is a sum of terms of
the form
\begin{equation*}
\bigl\{\{\underbrace{\_,\_,\_}_{k}\},\underbrace{\_,\_,\_,\_}_{l}\bigr\}
\end{equation*}
where $k+l=n$ and $k\gequ 1$. If there occur $4$ elements of $\Pi L$
or more, then among those $k$ or $l$ arguments there must be at
least $2$ in $\Pi L$, and it should be exactly $k=2$ and $l=2$,
since there are no brackets involving $3$ arguments in $\Pi L$. Then
the internal bracket also takes values in $\Pi L$, hence we get $3$
arguments in $\Pi L$ for the external bracket, so it must vanish.
Consider the Jacobiators that contain exactly $3$ arguments from
$\Pi L$. By a similar analysis one can see that the only potential
non-vanishing Jacobiator is  for $n=3$, which is exactly the
Jacobiator in $\Pi L$ and it vanishes since $L$ is a Lie
superalgebra. This leaves the Jacobiators with exactly $1$ or $2$
arguments in $\Pi L$. (The Jacobiators with all arguments in $V$
vanish by Theorem~\ref{thmjac} applied to $D$ such that $D^2=0$.)
They are as follows:
\begin{multline}\label{eqjacone1}
    J^{p+1}\left(\Pi x,a_1,\ldots,a_p\right)=\\
    \sum_{k=1}^p\sum_{\text{$(k,p-k)$-shuffles}}
    (-1)^{\xt+1}(-1)^{\e(\s;a_1,\ldots,a_p)}\,
    \left\{\Pi
    x,\{a_{\s(1)},\ldots,a_{\s(k)}\},a_{\s(k+1)},\ldots,a_{\s(p)}\right\}+\\
    \sum_{k=0}^p\sum_{\text{$(k,p-k)$-shuffles}}
    (-1)^{\e(\s;a_1,\ldots,a_p)}\,
    \left\{\{\Pi x,a_{\s(1)},\ldots,a_{\s(k)}\},a_{\s(k+1)},\ldots,a_{\s(p)}\right\}
\end{multline}
and
\begin{multline} \label{eqjactwo1}
    J^{p+2}(\Pi x,\Pi y,a_1,\ldots,a_p)=
    \left\{\{\Pi x,\Pi y\},a_1,\ldots,a_p\right\}+\\
    \sum_{k=0}^p\sum_{\text{$(k,p-k)$-shuffles}}\!\!\!(-1)^{\e(\s;a_1,\ldots,a_p)}
    \Biggl(
    (-1)^{\xt+1}\left\{\Pi x,\{\Pi y,a_{\s_(1)},\ldots,a_{\s(k)}\},
    a_{\s(k+1)},\ldots,a_{\s(p)}\right\} \\
    +(-1)^{(\yt+1)(\xt+\at_{\s(1)}+\ldots+\at_{\s(k)})}
    \left\{\Pi y, \{\Pi x,a_{\s_(1)},\ldots,a_{\s(k)}\},
    a_{\s(k+1)},\ldots,a_{\s(p)}\right\}
    \Biggr).
\end{multline}
Here $x,y\in L$, $a_i\in V$. By $(-1)^{\e(\s;a_1,\ldots,a_p)}$ we
denoted the sign arising from the action of a permutation $\s$ on
the product of $p$ commuting homogeneous variables of parities
$\at_1,\ldots,\at_p$. The equalities $J^{p+1}=0$ and $J^{p+2}=0$ can
be informally perceived, respectively,  as expressing the fact that
taking a bracket with $\Pi x$ acts, in a sense, as a derivation, and
that taking a bracket with $\{\Pi x, \Pi y\}$ acts, in a sense,  as
the commutator of brackets with $\Pi x$ and with $\Pi y$.  (All this
in a generalized sense, involving partitions and shuffles). Hence
these equalities are intuitively plausible. Let us prove them. For
this sake consider $x=y$ and $a_i=\x$ for all $i$, where $\xt=1$,
$\tilde\x=0$. Then~\eqref{eqjacone1} and \eqref{eqjactwo1} reduce to
\begin{multline}\label{eqjacone2}
    J^{p+1}(\Pi x,\x):=J^{p+1}(\Pi x,\underbrace{\x,\ldots,\x}_{p})=\\
    \sum_{k=1}^pC_p^k\bigl\{\Pi x, \{\underbrace{\x,\ldots,\x}_k\},
    \underbrace{\x,\ldots,\x}_{p-k}\bigr\}+
    \sum_{k=0}^pC_p^k\bigl\{\{\Pi x, \underbrace{\x,\ldots,\x}_k\},
    \underbrace{\x,\ldots,\x}_{p-k}\bigr\}
\end{multline}
and
\begin{multline}\label{eqjactwo2}
    J^{p+2}(\Pi x, \x):=J^{p+2}(\Pi x,\Pi x, \underbrace{\x,\ldots,\x}_{p})=\\
\bigl\{\{\Pi x,\Pi x\},\underbrace{\x,\ldots,\x}_p\bigr\}+
   2 \sum_{k=0}^p C_p^k
     \bigl\{\Pi x,\{\Pi x,\underbrace{\x,\ldots,\x}_k\},
    \underbrace{\x,\ldots,\x}_{p-k}\bigr\},
\end{multline}
respectively. Here $C_p^k$ denotes the binomial coefficient.
Substituting the definitions of the brackets in~\eqref{eqjacone2}, we
get after a simplification
\begin{multline*}
    J^{p+1}(\Pi x,\x)=\{-\Pi Dx+Px,\underbrace{\x,\ldots,\x}_{p}\,\}+\\
    \sum_{k=1}^pC_p^k \Bigl(\bigl\{\Pi x, P[\ldots[D\x,
    \underbrace{\x],\ldots,\x]}_{k-1},\underbrace{\x,\ldots,\x}_{p-k}\bigr\}
    +
    \bigl\{P[\ldots[[x,\underbrace{\x],\x],\ldots,\x]}_{k},
    \underbrace{\x,\ldots,\x}_{p-k}\bigr\}\Bigr)=\\
    (-1)^{p+1}P(\ad \x)^pDx +(-1)^pP(\ad\x)^pDPx+\\
\sum_{k=1}^pC_p^k \Bigl(
P(-1)^{p-1}(\ad\x)^{p-k}\left[x,P(\ad\x)^{k-1}D\x\right]
+P(-1)^{p}(\ad\x)^{p-k}DP(\ad\x)^kx \Bigr)
\end{multline*}
or
\begin{multline*}
    (-1)^pJ^{p+1}(\Pi x,\x)=
    -P(\ad \x)^pDx + P(\ad\x)^pDPx+\\
\sum_{k=1}^pC_p^k \Bigl(-
P(\ad\x)^{p-k}\left[x,P(\ad\x)^{k-1}D\x\right] +P
(\ad\x)^{p-k}DP(\ad\x)^kx \Bigr)\,.
\end{multline*}
Using the identity $(\ad\x)^kDP=-[(\ad\x)^{k-1}D\x,P(\,.\,)]$, for
$k\gequ 1$ (see the proof of Theorem~\ref{thmjac}), we can re-write
this as
\begin{multline}\label{eqjacone3}
    (-1)^pJ^{p+1}(\Pi x,\x)=P\Biggl(
    - (\ad \x)^pDx -[(\ad\x)^{p-1}D\x,Px]+\\
\sum_{k=1}^{p-1}C_p^k \Bigl(- \left[(\ad\x)^{p-k}
x,P(\ad\x)^{k-1}D\x\right] -
\left[(\ad\x)^{p-k-1}D\x,P(\ad\x)^kx\right] \Bigr)\\
-\left[x,P(\ad\x)^{p-1}D\x\right]+DP(\ad\x)^px\Biggr)=\\
    - P(\ad \x)^pDx +PD(\ad\x)^px-P[(\ad\x)^{p-1}D\x,x]-
\sum_{k=1}^{p-1}C_p^k P\left[(\ad\x)^{p-k}
x,(\ad\x)^{k-1}D\x\right]=\\
P\left[D,(\ad\x)^p\right]x-
\sum_{k=1}^{p}C_p^k P\left[(\ad\x)^{p-k} x,(\ad\x)^{k-1}D\x\right]
\end{multline}
where we used identities~\eqref{eqdistrib} and~\eqref{eqpdp}. Now, by
arguing in the same way as we did when deducing the
expression~\eqref{eqkomdiad1} for the commutator of $D$ and
$(\ad\x)^{N}$ acting on $D\x$ in the proof of Theorem~\ref{thmjac},
we can deduce the equality
\begin{multline*}
    \left[D,(\ad\x)^{p}\right]x=
\sum_{r=0}^{p-1}C_{p}^{p-1-r}\left[(\ad\x)^{r}D\x,
(\ad\x)^{p-1-r}x\right]=\\
\sum_{k=1}^{p}C_{p}^{k}\left[(\ad\x)^{k-1}D\x, (\ad\x)^{p-k}x\right].
\end{multline*}
Notice that since $x$ is odd, $\x$ is even, and $D$ is odd, in the
Lie bracket above both arguments are odd, so the order is irrelevant.
We immediately see that the two terms in the last line
of~\eqref{eqjacone3} cancel, and thus for all $x$  and $\x$
\begin{equation*}
J^{p+1}(\Pi x,\x)=0,
\end{equation*}
as desired. Now consider  $J^{p+2}(\Pi x,\x)$. Substituting the
definitions of the brackets into~\eqref{eqjactwo2}, we get
\begin{multline*}
    J^{p+2}(\Pi x,x)=
    -\bigl\{ \Pi [x,x],\underbrace{\x,\ldots,\x}_p\bigr\}+
   2 \sum_{k=0}^p C_p^k
     \bigl\{\Pi x,P(-\ad\x)^{k}x,
    \underbrace{\x,\ldots,\x}_{p-k}\bigr\}=\\
    -(-1)^pP(\ad\x)^p[x,x]+
    2(-1)^p\sum_{k=0}^pC_p^kP(\ad\x)^{p-k}\left[x,P(\ad\x)^kx\right],
\end{multline*}
or
\begin{multline*}
   (-1)^{p+1} J^{p+2}(\Pi x,x)=
     P(\ad\x)^p[x,x]-
    2 \sum_{k=0}^pC_p^kP(\ad\x)^{p-k}\left[x,P(\ad\x)^kx\right]=\\
 P(\ad\x)^p[x,x]-
    2 \sum_{k=0}^pC_p^kP\left[(\ad\x)^{p-k}x,P(\ad\x)^kx\right]=\\
    P(\ad\x)^p[x,x]-P\sum_{k=0}^pC_p^k\left[(\ad\x)^{p-k}x,(\ad\x)^kx\right]=\\
     P(\ad\x)^p[x,x]-P(\ad\x)^p[x,x]=0,
\end{multline*}
where we used the commutativity of $V$ and
identity~\eqref{eqdistrib}. Thus for all $x$  and $\x$
\begin{equation*}
J^{p+2}(\Pi x,\x)=0,
\end{equation*}
as desired. This completes the proof of the theorem.
\end{proof}

A remarkable fact about the formulae for the brackets in $\Pi
L\oplus V$ is that they arise naturally if one wants to extend  the
bracket in $\Pi L$ keeping the differential~\eqref{eqdifincoc2} a
derivation. Of course, the crucial  and much harder  thing is to
prove that they indeed give the structure of an $L_{\infty}$-algebra
as stated by Theorem~\ref{thmcocyl}. The subspace $V$ is a
subalgebra (even an ideal) with respect to this structure, and the
induced brackets are exactly the higher derived brackets.

\begin{corollary}\label{corcone}
The complex $\Pi L\oplus V$, with  operations defined as above, is a
cocylinder for $i\co \Pi K\to \Pi L$ in the category of
$L_{\infty}$-algebras. $V$ with the higher derived brackets of $D$
is a homotopy fiber (or a cocone), in this category, for the
inclusion of differential Lie superalgebras.
\end{corollary}

\begin{remark} The idea of relating the higher derived brackets of  $\D$
with homotopical algebra was proposed by the referee of the first
version of~\cite{tv:higherder}. He conjectured, for the $\ZZ$-graded
case, an interpretation of these brackets in terms of a `homotopy
left  cone' (cocone, in our terminology) and suggested a formula of
type~\eqref{eqbrackxaaa} for the extended brackets. In this section
we showed that the conjecture about a homotopical-algebraic
interpretation of  higher derived brackets is correct, in the
natural setup  where the brackets are generated by an arbitrary odd
derivation $D$. Corollary~\ref{corcone} gives the precise statement.
\end{remark}

The considerations of this section give an alternative and quite
unexpected, viewpoint of higher derived brackets. For a given
derivation $D$, which is assumed to be a differential, the
construction of the complex $\Pi L\oplus V$, viewed as a cone (for
$L\to V$) or a cocylinder (for $\Pi K\to \Pi L$) is canonical. The
higher derived brackets of $D$ appear as an answer to the question
of how to extend the algebra structure to $\Pi L\oplus V$ from $L$.

Notice also that although homological or homotopical algebra
requires $D^2=0$ from the start, we never directly used this
identity in the proof of Theorem~\ref{thmcocyl}, except where we
referred to Theorem~\ref{thmjac} in the particular case when
$D^2=0$; hence it seems reasonable that the homotopical-algebraic
picture can be rephrased in a way allowing to  incorporate a
possibly non-zero $D^2$.

\section{Generalizations and Discussion}

Let us return to Theorem~\ref{thmjac} and see what information can
be extracted from it if one does  not immediately set $D^2$ equal to
zero. To be able to make a precise statement, notice that our
construction of higher derived brackets allows  extension of
scalars, in the following sense.

Consider an arbitrary commutative superalgebra $\LL$ with unit (a
good example is the Grassmann algebra with $N$ generators,
$\LL=\LL_N$) and the tensor product $L\otimes \LL$. It is a Lie
superalgebra over $\LL$,  and we have $\Der_{\LL}(L\otimes
\LL)=(\Der L)\otimes \LL$. Thus the higher derived brackets can be
constructed from $D\in \Der_{\LL}(L\otimes \LL)$, i.e., a derivation
with coefficients in $\LL$. They will be operations on $V\otimes
\LL$. (In particular, brackets generated by $D\in \Der L$ can be
considered on $V\otimes \LL$ for any $\LL$ and this explains why it
is sufficient to check the Jacobiators only on even arguments.)
Clearly, Theorem~\ref{thmjac} remains valid. Now, the map which
assigns to a derivation $D$ all its higher derived brackets is a
linear operation in the sense that it commutes with sums and with
multiplication by scalars. Now we shall make use of the following
obvious algebraic statement: \textit{if a linear map of Lie
superalgebras maps the squares of odd elements to  squares, for all
extensions of scalars by various $\LL$, then it is a Lie algebra
homomorphism.} (Indeed, by polarization, it maps all brackets of odd
elements to the brackets; then by using suitable odd constants, even
elements can be turned into odd, and after that the constants can be
eliminated.)

An arbitrary sequence of multilinear symmetric operations on $V$ can
be encoded in a (formal) vector field $X$, which serves as their
generating function, so that the operations are obtained as the
higher derived brackets of $X$:
\begin{equation*}
\{u_1,\ldots,u_k\}_X=\left[\ldots\left[X,u_1\right],\ldots,u_k\right](0)
\end{equation*}
where $u_i\in V$, $X\in\Vect V$, as in~\eqref{equuQ}. If we restrict
ourselves to formal vector fields, this correspondence will be
one-to-one. The sequence of the Jacobiators of the brackets derived
from $X$ has the vector field $X^2$ as the generating function (this
is a very special case of Theorem~\ref{thmjac}, but can be seen
directly).

Consider now an arbitrary derivation $D\co L\to L$. Denote the vector
field on $V$ corresponding to the higher derived brackets of $D$, by
$Q_D$. Theorem~\ref{thmjac} then can be re-formulated as the equality
\begin{equation}\label{eqqdd}
(Q_D)^2=Q_{D^2}
\end{equation}
for all odd $D$. Having in mind the above remarks, we see that
Theorem~\ref{thmjac} is equivalent to the following.
\begin{theorem}
The correspondence $D\mapsto Q_D$ is a homomorphism of Lie
superalgebras $\Der L\to \Vect V$, i.e.,
\begin{equation}\label{eqqd1d2}
[Q_{D_1}, Q_{D_2}]=Q_{[D_1,D_2]}
\end{equation}
for all $D_1, D_2\in \Der L$.
\end{theorem}

(It is an interesting question whether there is a more direct way of
constructing a vector field on $V$ from the following data: the
homological field specifying the Lie bracket in $L$ and a derivation
$D$.)

Let $\mathfrak g$ be a Lie superalgebra and $V$ a vector space. We
call the space $V$ a \textit{generalized $L_{\infty}$-algebra over
$\mathfrak g$} (or: a \textit{$\mathfrak g$-parametric
$L_{\infty}$-algebra}) if there is given a homomorphism $\mathfrak
g\to \Vect V$. We can visualize this as (sequences of) brackets in
$V$ parametrized by elements of $\mathfrak g$. Relations between
elements of $\mathfrak g$ give rise to `generalized Jacobi
identities' in $V$ between the corresponding brackets.

\begin{example} If $\mathfrak g$ has dimension $0|1$, with a single
odd basis element $Q$ satisfying $Q^2=0$, then we get a usual
$L_{\infty}$-algebra structure.
\end{example}

\begin{example} If $\mathfrak g$ has dimension $1|1$, with a
basis $H,Q$ with $H$  even, $Q$  odd, satisfying $Q^2=H$, then a
generalized $L_{\infty}$-algebra over $\mathfrak g$ is the same as
an arbitrary sequence of odd symmetric brackets that a priori are
not subject to any relations. (In fact, there are some relations
that are always satisfied, they are the `mixed' Jacobi identities
for odd brackets and their Jacobiators, corresponding to the
identity $[H,Q]=0$.)
\end{example}

Apart from these two opposite extremes there should be other
interesting examples.

Another attractive direction is to  study  the higher derived
brackets where $V$ in the decomposition $L=K\oplus V$ is not assumed
Abelian. Notice that this is exactly the case in the original
definition of a (binary) derived bracket: given a Lie superalgebra
$L$ and an odd derivation $D\co L\to L$, then for arbitrary $a,b\in
L$
\begin{equation}\label{eqderived}
[a,b]_D:=[Da,b]
\end{equation}
(we use a sign convention convenient for the comparison
with~\eqref{eqdefderbrack}). This is a particular case
of~\eqref{eqdefderbrack} for $k=2$ if $L=V$ and $K=0$. It is known
that the derived bracket~\eqref{eqderived} is not, in general,
symmetric:
\begin{equation}\label{eqdersym}
[a,b]_D-(-1)^{\at\bt}[b,a]_D=D[a,b]
\end{equation}
(in typical applications it is possible to restrict  to an Abelian
subalgebra, thus restoring symmetry and making it into a different
special case of~\eqref{eqdefderbrack} for $k=2$ with a `hidden' $P$).

\begin{proposition} In the context of $L=K\oplus V$ where $V$ is
not necessarily Abelian, the $k$-th derived brackets defined
by~\eqref{eqdefderbrack} satisfy  the identity
\begin{multline}\label{eqddersym}
\{a_1,\ldots,a_i,a_{i+1},\ldots,a_k\}_D-(-1)^{\at_i\at_{i+1}}
\{a_1,\ldots,a_{i+1},a_i,\ldots,a_k\}_D=\\
\{a_1,\ldots,[a_i,a_{i+1}],\ldots,a_k\}_D
\end{multline}
for the transposition of two adjacent arguments, for all
$a_1,\ldots,a_k\in V$ and all $i=1,\ldots,k-1$. Here on the
right-hand side we have the $(k-1)$-th derived bracket with the Lie
bracket of the arguments $a_i$ and $a_{i+1}$ inserted at the $i$-th
position.
\end{proposition}

The proof is not hard and we omit it. Formula~\eqref{eqddersym}
generalizes~\eqref{eqdersym}.

It is known that the classical derived bracket, though not
symmetric, satisfies the Jacobi identity, defining  an odd Loday
algebra if $D^2=0$. What about analogs for higher derived brackets?
What is the precise list of relations in an algebraic structure
defined by the higher derived brackets if $V$ is non-Abelian? (It
includes an even Lie bracket as well as a sequence of odd brackets
and may be called an `$L_{\infty}$-algebra on a Lie algebra
background'.) It may be possible to make use of a
homotopic-algebraic approach such as in Section~\ref{sechomot}. We
hope to consider these questions elsewhere.

\appendix

\section{Appendix. Standard cylinders and cocylinders}
Here we collect, for reference purposes, the formulae for  the
standard constructions of   cylinders and cocylinders of chain maps
(compare, e.g., \cite{sgelman}). They all originate in topological
constructions of the cylinder $X\times I$ and cocylinder $X^I$.

A \textit{complex} is a ($\Z$-graded) vector space equipped with an
odd operator $d$ such that $d^2=0$. A \textit{map} or a `chain map'
is an even linear map commuting with $d$.

Let $f\co X\to Y$ be a map of complexes.

The standard \textit{cylinder} diagram for $f\co X\to Y$ is the
commutative diagram
\begin{equation*}
\begin{diagram}[small]
    X &       & \rTo^f          &       & Y  \\
      & \rdTo_j &               & \ruTo_p &    \\
      &       & \Cyl f &       &    \\
\end{diagram}
\end{equation*}
where
\begin{equation*}
\Cyl f=X\oplus \Pi X \oplus Y
\end{equation*}
with the differential given by
\begin{equation*}
d(x_1,x_2,\Pi y)=(dx_1-x_2, \Pi (-dx_2), dy+f(x_2)).
\end{equation*}
The maps $j$ and $p$ are given by the formulae
\begin{align*}
    j(x)&=(x,0,0)\\
    p(x_1,\Pi x_2,y)&=f(x_1)+y,
\end{align*}
and $p$ is a quasi-isomorphism with a quasi-inverse map $i\co Y\to
\Cyl f$, $i(y)=(0,0,y)$. The \textit{cone} of $f$ is the cofiber of
$j$, i.e., $\Cyl f/j(X)$. Hence
\begin{equation*}
\Cone f=\Pi X\oplus Y
\end{equation*}
with the differential
$$
d(\Pi x,y)=(\Pi (-dx), dy+f(x)).
$$

In a similar way, the standard \textit{cocylinder} diagram for $f\co
X\to Y$ is the commutative diagram
\begin{equation*}
\begin{diagram}[small]
    X &       & \rTo^f          &       & Y  \\
      & \rdTo_j &               & \ruTo_p &    \\
      &       & \Cocyl f &       &    \\
\end{diagram}
\end{equation*}
where
\begin{equation*}
\Cocyl f=X\oplus   Y \oplus \Pi Y
\end{equation*}
with the differential given by
\begin{equation*}
d(x,y_1,\Pi y_2)=(dx, dy_1, \Pi (f(x)-y_1-dy_2)).
\end{equation*}
The maps $j$ and $p$ are given by the formulae
\begin{align*}
    j(x)&=(x,f(x),0)\\
    p(x,y_1,\Pi y_2)&=y_1,
\end{align*}
and $j$ is a quasi-isomorphism with a quasi-inverse map $q\co \Cocyl
f\to X$, $q(x,y_1,\Pi y_2)=x$. The \textit{cocone} of $f$ is the
fiber (kernel) of $p$. Hence
\begin{equation*}
\Cocone f=X\oplus \Pi  Y
\end{equation*}
with the differential
\begin{equation*}
d(x,\Pi  y)=(dx, \Pi (f(x)-dy)).
\end{equation*}

It follows that $\Pi \Cone f=\Cone f^{\Pi}=\Cocone (-f)$; i.e., up to
a sign, the cone and cocone of a chain map $f$ are related by the
parity shift functor. In the main text, the complex $L\oplus \Pi V$
appearing there as a cocylinder of the inclusion of complexes $i\co
K\to L$, can be alternatively viewed  as the canonical $\Cocone (-P)$
or as $\Pi \Cone P$ where the projector $P$ is treated as a map $L\to
V$, so $V$ with the differential $PD$ is considered as a quotient
complex of $L$ (rather than a subspace of $L$).




\def\cprime{$'$} \def\cprime{$'$}

\noindent
      Dr Theodore Voronov \\
      School of Mathematics\\University of Manchester\\Sackville
      Street\\
Manchester M60 1QD\\United Kingdom\\
      \texttt{theodore.voronov@manchester.ac.uk}

\label{lastpage}
\end{document}